\documentclass{amsart}
\usepackage{enumerate}
\usepackage{amsfonts}
\usepackage{amsmath}
\usepackage{amssymb}
\usepackage{amsthm}
\usepackage{url}

\newtheorem{teor}{Theorem}[section]
\newtheorem{cor}[teor]{Corollary}

\newtheorem{rem}[teor]{Remark}
\newtheorem{lem}[teor]{Lemma}
\newtheorem{prop}[teor]{Proposition}
\newtheorem{defi}[teor]{Definition}

\newcommand{\Z}{\mathbb Z}

\newcommand{\calA}{\mathcal{A}}

\def\Jac{\textrm{Jac}}

\def\op{\operatorname}

\def\as{\op{AS}}
\def\ars#1{\renewcommand\arraystretch{#1}}
\def\aut{\op{Aut}}

\def\dd{{\mathcal D}}

\def\dm{\op{dim}}

\def\diso{\lower.4ex\hbox{$\downarrow$}\raise.4ex\hbox{\mbox{\scriptsize $\wr$}}}

\def\enc{E_{\op{nc}}}
\def\endo{\op{End}}

\def\eps{\epsilon}

\def\ff#1{\mathbb F_{q^{#1}}}

\def\fq{\mathbb F_q}
\def\ft{\mathbb F_2}
\def\ff{\mathbb F_4}
\def\fe{\mathbb F_8}

\def\gen#1{\big\langle\, {#1} \,\big\rangle}

\def\hm#1#2#3{\op{Hom}_{#1}(#2,#3)}

\def\k{\op{Ker}}
\def\kb{\bar{k}}
\def\kc{(k^*)^3}
\def\la{\lambda}
\def\lg{l\raise.6ex\hbox to.2em{\hss.\hss}l}
\def\lra{\longrightarrow}

\def\md#1{\ \,\mbox{\rm(mod }{#1})}

\def\orb{\hbox to  .3em{$\backslash$}\backslash}

\def\tq{\,\,|\,\,}
\def\tr{\op{Tr}}

\newcounter{cs}
\stepcounter{cs}
\newcommand{\casos}{\begin{itemize}}
\newcommand{\fcasos}{\end{itemize}\setcounter{cs}{1}}

\newfont{\tit}{cmr12 scaled \magstep3}

\title[Jacobians in isogeny classes]{Jacobians in isogeny classes of supersingular abelian threefolds in characteristic $2$}

\author{Enric Nart}\thanks{The first author acknowledges support from the project MTM2006-11391 from the Spanish MEC}
\address{Departament de Matem\`atiques,
        Universitat Aut\`onoma de Barcelona,
        Edifici C, 08193 Bellaterra, Barcelona, Spain.}
\email{nart@mat.uab.cat}

\author{Christophe Ritzenthaler}
\address{Institut de Math\'ematiques de Luminy,
         UMR 6206 du CNRS,
         Luminy, Case 907, 13288 Marseille, France.}
\email{ritzenth@iml.univ-mrs.fr}

\date{14 September 2006}



\begin{document}

\maketitle

\begin{abstract}
We exhibit  the isogeny classes of supersingular abelian threefolds over
$\mathbb{F}_{2^n}$ containing the Jacobian of a genus $3$ curve. In
particular, we prove that for even $n>6$ there always exist a
maximal and a minimal curve over $\mathbb{F}_{2^n}$. All the curves
can be obtained explicitly.
\end{abstract}

\section*{Introduction}
Let $g \geq 0$ be an integer and $k=\fq$ a finite field with
$q=p^n$. When $C$ spans the finite set of (smooth, absolutely
irreducible, projective) curves of genus $g$ over $\fq$ what is its
maximal number of points ? This question has an easy answer in the
case $g=0$ and $1$ but one has to wait till 1985 for Serre to solve
the case $g=2$ \cite{serre_point}. For $g \geq 3$, the problem has
only partial solutions, most of the time obtained case by case for
small $q$ \cite{vdgeer}. This problem can be seen as part of a more
general question : let $\calA$ be an isogeny class of abelian
varieties of dimension $g \geq 1$ over $k$. Does $\calA$ contain the
Jacobian of a curve of genus $g$ over $k$ ? Again, the answer to
this question, easy in the case $g=1$, has been obtained
only recently for $g=2$ \cite{hnr}. For $g \geq 4$, the problem seems more
 difficult even in the geometric setting (i.e. over $\overline{k}$)
since one has to
characterize Jacobians among abelian varieties.\\
So, let us concentrate on the new frontier, the case $g=3$. As for
$g=2$, if
$(A,\lambda)$ is an undecomposable principally polarized abelian
threefold defined over $k$ then there exists a curve $C$ defined
over $k$ such that $\Jac(C)$ is isomorphic to $(A,\lambda)$ over
$\overline{k}$ \cite{oort23}. However, over $k$ the situation is more
subtle than it was for $g=2$ \cite[Appendix]{lauter} :
\begin{itemize}
\item if $C$ is hyperelliptic then $\Jac(C) \simeq A$.
\item if $C$ is non-hyperelliptic then $\Jac(C)$ is either
isomorphic to $A$ or to the quadratic twist of $A$.
\end{itemize}
Thus, if $A$ belongs to the isogeny class with Frobenius polynomial
$x^6+a x^5+ b x^4+ c x^3+ q b x^2+ a  q^2 x+ q^3$, which we denote
$\calA_{(a,b,c)}$, then $\Jac(C)$ may belong to this class or to the
class $\calA_{(-a,b,-c)}$. In particular $\# C(k)=1+q \pm a$. So
far, nobody knows a practical method to determine which case occurs.
In account of this uncertainty on the sign, the current best general result  is
due to Lauter \cite{lauterg3} : let $m=\lfloor 2 \sqrt{q} \rfloor$
then there exists a curve $C/k$ of genus $3$ such that
$$\#C(k) \geq q+1+3m -3 \quad \textrm{or} \quad \#C(k) \leq q+1-3m
+3.$$ If we restrict to the case $n$ even, Ibukiyama
\cite{ibukiyama} showed that if $p \equiv 3 \pmod{4}$ or ($n \equiv
2 \pmod{4}$ and $p \neq 2$) then there exists a maximal curve over
$k$ (recall that a curve is {\em maximal} (resp. {\em minimal})
if $\#C(k)=q+1+g m$ (resp. $\#C(k)=q+1-g m$)). \\
In this article we consider the case $p=2$ and we answer the
following question: which are the isogeny classes of {\em
supersingular} abelian threefolds over $k=\mathbb{F}_{2^n}$
containing the Jacobian of a curve ? Theorems \ref{thsimple},
\ref{thhalf}, \ref{teorfinal} give a complete answer, which is easy
to state for $n>6$:\bigskip

{\bf Result.} {\em All isogeny classes of supersingular abelian
threefolds over $k$ with Frobenius polynomial $F$ contain the
Jacobian of a genus $3$ curve except in the following cases (where
$\eps=\pm 1$ takes both values):
\begin{itemize}
\item $n$ even and $F(x)=f(x) (x^4+\eps\sqrt{q}x^3+qx^2+\eps q\sqrt{q}x+q^2)$
with $f$ the Frobenius polynomial of a supersingular elliptic curve
over $k$ ;
\item $n$ odd and $F$ is equal to
$(x^2+q)(x^4+\eps\sqrt{2q}x^3+qx^2+\eps q\sqrt{2q}x+q^2)$ or
$(x^2+\eps\sqrt{2q}x+q)(x^4+\eps\sqrt{2q}x^3+qx^2+\eps
q\sqrt{2q}x+q^2)$.
\end{itemize}}\bigskip

In particular, the above mentioned ambiguity of sign disappears when
we consider supersingular abelian threefolds up to isogeny and when
the base field is big enough ($n>6)$. In this case, the isogeny
class $\calA_{(a,b,c)}$ contains a Jacobian if and only if the
isogeny class $\calA_{(-a,b,-c)}$ does. One may speculate if a
similar phenomenon will occur for non-supersingular threefolds
and/or other characteristics.

 As another by-product, we obtain the following result.\bigskip

{\bf Corollary.} {\em  For $n>6$ even,  there are maximal curves and
minimal curves of genus three over $k$. For $n=6$ there are maximal
curves over $k$, but no minimal curve. For $n=2$ or $4$
there are neither maximal nor minimal curves over $k$.}\bigskip

Let us point out that it is easy to use the constructive methods of
the proofs to build explicit models for these curves. Indeed, our
 approach is based on  models for supersingular quartics $C$
over $k$ (genus $3$ non-hyperelliptic curves whose Jacobians are
supersingular)  \cite{nr} :
\begin{equation*}
C\colon\quad Y^4+fY^2+gY=X^3+dX^2+e, \qquad g\in k^*,\ d,\,e,\,f\in
k.
\end{equation*}
 The automorphism
group of $C/\overline{k}$ contains $G=\Z/2\Z \times \Z/2\Z$, and the
quotients of $C$ by the three non trivial automorphisms of $G$
define three explicit elliptic curves $E_i/\overline{k}$ such that
$\Jac(C) \sim_{\overline{k}} E_1 \times E_2 \times E_3$. We reverse
then the process and starting with a product $E_1 \times E_2 \times
E_3$ of supersingular elliptic curves over $k$ (or a quadratic or
cubic extension of $k$), we determine if it is isogenous to the
Jacobian of a curve of $C$ of this type with
coefficients in $k$.

This paper is organized in the following way : in Sec.~1, we present
elementary results on supersingular elliptic curves over $k$,
expressing their trace in terms of the coefficients of normalized
models. In Sec.~2 we give parametrizations of some genus $0$ curves
and study the number of points of certain hyperelliptic curves of
genus $4$ as well as of an affine surface. These results are
useful in Sec.~3, which is the core of the paper. In  Sec.~3 we
determine existence of supersingular quartics with elliptic quotients in prescribed isogeny classes.
This section is divided into three parts according to
the splitting behavior of $Y^3+f Y+g$. Finally, Sec.~4 answers our
initial question by comparing the results obtained in Sec.~3 and the
isogeny classes of supersingular abelian threefolds over $k$.

\section{Supersingular elliptic curves in characteristic
$2$.}\label{ssec} Let $k=\fq$ be a finite field of characteristic
$2$, with $q=2^n$.
We shall denote simply by $\tr$ or $\tr_k$ the absolute trace
$\tr_{k/\ft}$.

The perfect duality of $k$ given by the nondegenerate pairing
$$
\gen{\ ,\ }\colon k\times k \lra \ft,\qquad \gen{x,y}=\tr(xy),
$$
leads to a canonical isomorphism of $\ft$-vector spaces
$$
(\ \,)^*\colon k\lra \hm{\ft}{k}{\ft}, \qquad x\mapsto x^*=\tr(x-).
$$
Moreover, if $q$ is a square we have an analogous result
substituting $\ft$ by $\ff$. We shall often use the following
immediate consequence of this duality
\begin{lem}\label{onto}
Let $k_0=\ft$ or $k_0=\ff$ (if $q$ is a square). If $u_1,\dots,u_r\in k$ are $k_0$-linearly independent, then we have an onto map:
$$
k\lra (k_0)^r,\qquad x\mapsto (x^*(u_1),\dots,x^*(u_r)).
$$
\end{lem}

We recall some basic facts concerning the Artin-Schreier
operators:

$$ \as(x)=x+x^2,\qquad \as^2(x)=x+x^4.
$$
There is an exact sequence of additive abelian
groups $$ 0\lra \ft\stackrel{i}\lra k\stackrel{\as}\lra
k\stackrel{\tr}\lra\ft\lra 0,
$$
and a similar exact sequence of $\ff$-vector spaces, when
$\ff\subseteq k$
$$
0\lra \ff\stackrel{i}\lra k\stackrel{\as^2}\lra
k\stackrel{\tr_{k/\ff}}\lra\ff\lra 0.
$$

\begin{lem}\label{qeq}
A quadratic polynomial $f(x)=x^2+ax+b\in k[x]$ is separable if and
only if $a\ne0$; in this case, $f(x)$ is irreducible if and only if
$b/a^2\not\in\as(k)$.
\end{lem}

\subsection*{Artin-Schreier models}
An elliptic curve over $k$ is supersingular if and only if its $j$
invariant equals $0$. Thus, the supersingular elliptic curves are
all geometrically isomorphic and they have Weierstrass models:
$$
y^2+a_3y=x^3+a_2x^2+a_4x+a_6, \quad a_3\ne0.
$$

The map $(x,y)\mapsto (x,a_3^{-1}y)$ sets a $k$-isomorphism with
an Artin-Schreier model
$$
y^2+y=a_3^{-2}x^3+a_3^{-2}a_2x^2+a_3^{-2}a_4x+a_3^{-2}a_6.
$$
A change of variables $y=y+a_3^{-2}a_4x+t$, $t\in k$, makes the coefficient of
$x$ vanish and allows us to suppose that the last term is either $0$ or
a fixed element $c_0$ in $k\setminus \as(k)$. We are led in this way to normal Artin-Schreier models of the type
$$
E\colon\quad y^2+y=ax^3+bx^2+c, \quad a\in k^*,\ b\in k,\ c\in\{0,c_0\}.
$$

The advantage of the Artin-Schreier models is that the canonical
twist (associated to the automorphism ``multiplication by $-1$")
is obtained just by making a different choice of the coefficient $c$.

$$
E'\colon\quad y^2+y=ax^3+bx^2+c+c_0.
$$
If $x^2+tx+q$ is the Frobenius polynomial of $E/k$ then $x^2-tx+q$
is the Frobenius polynomial of $E'/k$.

 Let us introduce a
special notation for some particular curves
$$\ars{1.4}
\begin{array}{cll}
E_a\colon &\quad y^2+y\,=\,ax^3,  &\quad  a\in k^*,\\
H\colon &\quad y^2+y\,=\,x^3+x^2,& \\
E_0\colon &\quad y^2+y\,=\,c_0^{-3}x^3+x^2,&\quad c_0\in k\setminus \as(k).
\end{array}\ars{1}
$$

\begin{teor}\label{allcurves}
There are seven or three $k$-isomorphism classes of supersingular
elliptic curves according to $q$ being a square or not. In Tables
1 and 2 we find representative models of these classes, and their
respective Frobenius polynomials $x^2+tx+q$.
\end{teor}

\begin{table}
\begin{center}
\renewcommand{\arraystretch}{1.4}
\begin{tabular}{|c|c|c|c|c|c|}
\hline $E$&$E_1$&$E_1'$&$E_a$, $E_{a^2}$&$E_a'$, $E_{a^2}'$&$E_0$\\
\hline$t$&$-(-1)^{n/2}2\sqrt{q}$&$(-1)^{n/2}2\sqrt{q}$&$(-1)^{n/2}\sqrt{q}$&$-(-1)^{n/2}\sqrt{q}$&0
\\\hline\end{tabular}\vskip.2cm
\caption{Supersingular
elliptic curves, $q$ square, $a\in k^*\setminus \kc$}
\label{Table:EllSq}
\end{center}
\end{table}

\begin{table}
\begin{center}\ars{1.4}
\begin{tabular}{|c|c|c|c|}
\hline
$E$&$E_1$&$H$&$H'$\\
\hline$t$&$0$&$(-1)^{(n^2-1)/8}\sqrt{2q}$&$-(-1)^{(n^2-1)/8}\sqrt{2q}$
\\
\hline\end{tabular}\vskip.2cm
\caption{Supersingular
elliptic curves, $q$ non-square}
\label{Table:EllNsq}
\end{center}
\end{table}

This can be checked by computing all
$\kb/k$ twists of the curve $\,E_1\colon\ y^2+y\,=\,x^3$. If a twist $E$
corresponds to an automorphism $\alpha\in\aut(E_1)$ and $f\colon
E_1\lra E$ is a geometric isomorphism with
$\alpha=f^{-1}f^{\sigma}$, then the pullback by $f$ of the
Frobenius endomorphism of $E$ is $\alpha\pi$, where $\pi$ is the
Frobenius endomorphism of $E_1$. Since $\pi^2=(-1)^nq$, this determines immediately the
Frobenius polynomial of $E$ in terms of the characteristic
polynomial of $\alpha$ as an element of $\endo(E)$.


\begin{prop}\label{crit}
Let $E$ be a supersingular elliptic curve over $k$ given by an Artin-Schreier
normal model
$$E\colon\quad y^2+y=ax^3+bx^2+c, \quad a\in k^*,\
b\in k,\ c\in \{0,c_0\}.
$$

Tables 3 and 4 below determine the isomorphism class of $E$ in
terms of the parameters $a$, $b$ and $c$.
\end{prop}

\begin{table}
\begin{center}\ars{1.4}
\begin{tabular}{|c|c|}
\hline $a$, $b$&$E$ isomorphic to \\
\hline $\begin{array}{c}a\not\in\kc\mbox{ or}\\
a=u^3,\
bu^{-2}\in\as^2(k)\end{array}$&$\begin{array}{ll}E_a,&\mbox{ if
}\tr(c+v^3a)=0\\E_a',
&\mbox{ if }\tr(c+v^3a)=1\end{array}$\\
\hline $a=u^3,\ bu^{-2}\not\in\as^2(k)$&$E_0$\\\hline
\end{tabular}\vskip.2cm
\caption{Isomorphism class of $E\colon y^2+y=ax^3+bx^2+c$,\qquad\qquad\quad
$q$ square, \  $v$ is any solution in $k$ to \ $va+v^4a^2=b$}
\label{Table:CritSq}
\end{center}
\end{table}

\begin{table}
\begin{center}\ars{1.4}
\begin{tabular}{|c|c|}
\hline $a$, $b$&$E$ isomorphic to \\
\hline$a=u^3,\
bu^{-2}\in\as(k)$&$E_1$\\
\hline$ a=u^3,\ bu^{-2}\not\in\as(k)$&$\begin{array}{ll}\
H,&\mbox{ if }\tr(c+v^3+v)=0\\H',
&\mbox{ if }\tr(c+v^3+v)=1\end{array}$\\
\hline
\end{tabular}\vskip.2cm
\caption{Isomorphism class of $E\colon y^2+y=ax^3+bx^2+c$,\qquad\qquad\quad
$q$ non-square, \
  $v$ is any solution in $k$ to \ $1+v+v^4=bu^{-2}$}
\label{Table:CritNsq}
\end{center}
\end{table}

\begin{proof}
A $k$-isomorphism from $E$ to another normal model $\tilde{E}$ with
parameters $\tilde{a},\,\tilde{b},\,\tilde{c}$ is necessarily of the form:
$$
(x,y)\mapsto (u(x+v),y+av^2x+t),
$$
with $u\in k^*$, $v,\,t\in k$ satisfying
$$
\tilde{a}=u^{-3}a,\quad \tilde{b}=u^{-2}(b+a^2v^4+av),\quad \tilde{c}=
c+a^2v^6+bv^2+t+t^2.
$$
In particular, the class of $a\in k^*/\kc$ is preserved.

We can get a target curve with $\tilde{b}=0$ if and only if $a^2v^4+av=b$
for some $v\in k$. Note that in this case $bv^2=a^2v^6+av^3$
belongs to $\as(k)$ and $\tilde{c}\equiv c+av^3\md{\as(k)}$. If $a\not\in
\kc$ the $\ft$-linear map $x\mapsto a^2x^4+ax$ is an automorphism
of $k$ and the equation $a^2v^4+av=b$ has a unique solution $v\in
k$. If $a=z^3$ for some $z\in k$, this equation is equivalent to
$w^4+w=bz^{-2}$ for $w=zv$, and it has a solution in $k$ if and
only if $bz^{-2}\in\as^2(k)$.

We can get a target curve with $\tilde{a}=\tilde{b}=1$ if and only if $a=u^3$
and $a^2v^4+av=b+u^2$ for some $v\in k$, or equivalently (for
$w=uv$): $w^4+w=bu^{-2}+1$. In this case $av^3=w^3$ and
$bv^2=a^2v^6+av^3+w^2$, so that $\tilde{c}\equiv c+w^3+w\md{\as(k)}$.

After these considerations, the proposition is consequence of Theorem \ref{allcurves}.
\end{proof}

\section{Some special curves and surfaces}\label{scs}
\begin{lem}\label{parameterization}
\begin{enumerate}[a)]
\item
The affine plane curve
$$
a+bx+c^4x^4=a'+b'y+(c')^4y^4, \qquad b'c\ne bc',
$$
has genus zero and it admits the parametrization:
$$
x=\dfrac{c'(a+a'+\mu^4)+b'\mu}{b'c+bc'},\qquad y=\dfrac{c(a+a'+\mu^4)+b\mu}{b'c+bc'}
$$
in terms of a free parameter $\mu$.
\item
For $r,s,t\in \kb^*$ such that $r+s+t=0$, the affine spatial curve
$$
x+x^4r^8=y+y^4s^8=z+z^4t^8
$$
has genus zero and it admits the parametrization:
\begin{equation}\label{param1}
\ars{1.4}
\begin{array}{l}
x=(rst)^{-2}\left(rst\la+(r^{-2}t^4+rs)\la^4+r^{-2}\la^{16}\right),\\
y=(rst)^{-2}\left(rst\la+(s^{-2}r^4+st)\la^4+s^{-2}\la^{16}\right),\\
z=(rst)^{-2}\left(rst\la+(t^{-2}s^4+tr)\la^4+t^{-2}\la^{16}\right),
\end{array}
\end{equation}
in terms of a free parameter $\la$. In particular, for any point $(x,y,z)$ of this curve
\begin{multline*}
x^3r^8+y^3s^8+z^3t^8=(rst)^{-4}\la^{36}+(rst)^{-2}\la^{18}+\\+(rst)^{-4}(r^7s+s^7t+t^7r)\la^{12}+
(rst)^{-1}\la^9.
\end{multline*}
\item
For $r,s,t\in \kb^*$ such that $r^{-3}+s^{-3}+t^{-3}=0$, the affine spatial curve
$$
r^{16}(1+x+x^4)=s^{16}(1+y+y^4)=t^{16}(1+z+z^4)
$$
has genus zero and it admits the parametrization:
\begin{equation}\label{param2}\ars{1.4}
\begin{array}{l}
x=r^{-4}\left[(r^{-1}+s^{-1}+t^{-1})^8(r^{12}+s^{12}+t^{12}+(st)^6)+\right.\\
\qquad\qquad \quad\left.+r^{-9}(st)^3\la+(r^{-18}s^{-6}+t^{-24})\la^4+(rst)^{-36}\la^{16}\right],\\
y=s^{-4}\left[(r^{-1}+s^{-1}+t^{-1})^8(r^{12}+s^{12}+t^{12}+(tr)^6)+\right.\\
\qquad\qquad \quad \left.+s^{-9}(tr)^3\la+(s^{-18}t^{-6}+r^{-24})\la^4+(rst)^{-36}\la^{16}\right],\\
z=t^{-4}\left[(r^{-1}+s^{-1}+t^{-1})^8(r^{12}+s^{12}+t^{12}+(rs)^6)+\right.\\
\qquad\qquad \quad \left.+t^{-9}(rs)^3\la+(t^{-18}r^{-6}+s^{-24})\la^4+(rst)^{-36}\la^{16}\right],
\end{array}
\end{equation}
in terms of a free parameter $\la$. In particular, for any point $(x,y,z)$ of this curve
\begin{multline*}
x^3+x+y^3+y+z^3+z=R^{-84}\la^{36}+R^{-72}S^8\la^{32}+R^{-42}\la^{18}+
R^{-36}S^4\la^{16}+\\+R^{-12}T\la^{12}+R^{-21}\la^9+S^8T\la^8+
R^{12}S^{64}\la^4+R^{24}S^{72},
\end{multline*}
where $R=rst$, $S=r^{-1}+s^{-1}+t^{-1}$, and $T=r^{-6}s^{-42}+s^{-6}t^{-42}+t^{-6}r^{-42}$.
\end{enumerate}
\end{lem}

\begin{proof}
The only singularity of the curve $a+bx+c^4x^4=a'+b'y+(c')^4y^4$ is a triple point at infinity. Thus, it can be parametrized in terms of the pencil of lines through this singular point. This explains the first item.

In particular, we can express the solutions of $x+x^4r^8=y+y^4s^8$, and independently $x+x^4r^8=z+z^4t^8$, in terms of free parameters
$$
\ars{1.6}
\begin{array}{l}
x=t^{-2}(\mu+\mu^4s^2),\quad y=t^{-2}(\mu+\mu^4r^2), \\
x=s^{-2}(\nu+\nu^4t^2),\quad z=s^{-2}(\nu+\nu^4r^2).
\end{array}
$$
The solutions of  $x+x^4r^8=y+y^4s^8=z+z^4t^8$ correspond to $t^{-2}(\mu+\mu^4s^2)=x=s^{-2}(\nu+\nu^4t^2)$.  By the first item we can express $\mu$, $\nu$ in terms of a free parameter $\rho$
$$
\mu=\dfrac{s^{-2}(st)^{1/2}\rho+t\rho^4}{s^{-1}+t^{-1}},\qquad \nu=\dfrac{t^{-2}(st)^{1/2}\rho+s\rho^4}{s^{-1}+t^{-1}}.
$$
The change of parameters $\la=(st)^{1/2}\rho$ leads to the desired expression of $x,\,y,\,z$ in terms of $\la$. The computation of $x^3r^8+y^3s^8+z^3t^8$ in terms of $\la$ is straightforward.

In a similar way, we can express the solutions of $r^{16}(1+x+x^4)=s^{16}(1+y+y^4)$, and independently $r^{16}(1+x+x^4)=t^{16}(1+z+z^4)$, in terms of free parameters
$$
\ars{2.4}
\begin{array}{l}
x=r^{-4}\dfrac{r^{16}+s^{16}+s^{12}\mu+\mu^4}{r^{12}+s^{12}},\quad y=s^{-4}\dfrac{r^{16}+s^{16}+r^{12}\mu+\mu^4}{r^{12}+s^{12}}, \\
x=r^{-4}\dfrac{r^{16}+t^{16}+t^{12}\nu+\nu^4}{r^{12}+t^{12}},\quad
z=t^{-4}\dfrac{r^{16}+t^{16}+r^{12}\nu+\nu^4}{r^{12}+t^{12}}.
\end{array}
$$
The solutions of  $r^{16}(1+x+x^4)=s^{16}(1+y+y^4)=t^{16}(1+z+z^4)$ correspond to values of $\mu,\,\nu$ such that the two above expressions for $x$ coincide. Using the identity $(rs)^3+(rt)^3+(st)^3=0$ this is equivalent to
\begin{multline*}
(rt)^{12}(r^4+s^4)+(rt)^{12}\mu+(r^{12}+t^{12})\mu^4=\\=(rs)^{12}(r^4+t^4)+(rs)^{12}\nu+(r^{12}+s^{12})\nu^4.
\end{multline*} By the first item we can express $\mu$, $\nu$ in terms of a free parameter $\la$
$$\ars{1.4}
\begin{array}{l}
\mu=t^{-6}\left[(rst)^6(r^{-1}+s^{-1}+t^{-1})^8+r^3s^3t^{-3}\la+(rst)^{-6}\la^4\right],\\
\nu=s^{-6}\left[(rst)^6(r^{-1}+s^{-1}+t^{-1})^8+r^3t^3s^{-3}\la+(rst)^{-6}\la^4\right].
\end{array}
$$
The computation of $x,\,y,\,z$ and $x^3+x+y^3+y+z^3+z$ in terms of $\la$ is then straightforward.
\end{proof}

We need to control the number of points of supersingular hyperelliptic curves of genus $4$ given by equations of the type
$$
\dd_{ABC}\colon\quad y^2+y=Ax^9+Bx^3+Cx,\qquad A\in k^*,\,B,\,C\in k.
$$

\begin{lem}\label{dabc}
Consider the $\ft$-linear endomorphism of $k$ given by the polynomial $P(x)=A^8x^{64}+B^8x^{16}+B^4x^4+Ax\in k[x]$. Let $W=\k(P)\cap k$, $w=\dm_{\ft}(W)$ and
let $Q\colon W\lra \ft$ be the linear map given by $Q(x)=\tr(Ax^9+Bx^3+Cx)$.
Then,
$$\ars{1.2}
|\dd_{ABC}(k)|=\left\{\begin{array}{ll}
q+1\pm \sqrt{2^wq},&\quad\mbox{ if }Q=0,\\
q+1,&\quad\mbox{ if }Q\ne0.
\end{array}\right.
$$
In particular, if $w<n$ we have $2\le |\dd_{ABC}(k)|\le 2q$. Moreover, one has $w=n$ only in the following cases
\begin{itemize}
\item $q=64$, $B=0$, $A\in\fe$,
\item $q=16$, $B=A^2$,
\item $q=8$, $B=0$,
\item $q=4$, $A+B\in\ft$,
\item $q=2$.
\end{itemize}
Finally, if $q$ is a square, $C=0$ and $w=n$ we have $|\dd_{ABC}(k)|=2q+1$.
\end{lem}

\begin{proof}
The statement about $|\dd_{ABC}(k)|$ is consequence of \cite[secs. 3,5]{vv}. The condition $w=n$ means that the polynomial $P(x)$ vanishes identically on $k$ and it is easy to check that this translates into the conditions given in the lemma. The last statement is consequence of the fact that the linear form $\tr(Ax^9+Bx^3)$ vanishes identically if $w=n$; for instance, if $q=16$:
$$
Ax^9+A^2x^3\equiv Ax^9+(A^2x^3)^8=Ax^9+Ax^9=0\md{\as(k)}.
$$
\end{proof}

\subsection*{Number of triples of cubes adding up to zero}
\begin{lem}\label{surface}
The affine surface $x^3+y^3+z^3=0$ has $(q-1)\left(|E_1(k)|-3\mu\right)$ $k$-rational points with nonzero coordinates, where $\mu=|\mu_3(k)|$ is equal to $3$ or $1$ according to $q$ being a square or not.
\end{lem}

\begin{proof}
The plane projective curve given in homogeneous coordinates by $x^3+y^3+z^3=0$ is $\ft$-isomorphic to $E_1$ through the map $(x,y,z)\mapsto (x+y+z,x+z,x+y)$. This Fermat curve has $3\mu$ points satisfying $xyz=0$, and each projective solution with $xyz\ne0$ determines $q-1$ affine points of the surface $x^3+y^3+z^3=0$.
\end{proof}

\section{Supersingular plane quartics with prescribed elliptic quotients}
In \cite{nr} it is proved that any supersingular plane quartic defined over $k$ admits an affine model of the type
\begin{equation}\label{quartic}
C\colon\quad Y^4+fY^2+gY=X^3+dX^2+e, \qquad g\in k^*,\ d,\,e,\,f\in k,
\end{equation}
with a single point $P_{\infty}$ at infinity, which is a hyperflex.
The line at infinity is the only bitangent of $C$. For each nonzero
root $\theta$ of the polynomial $R(Y):=Y^4+fY^2+gY$, the quotient of
$C$ by the involution $Y\mapsto Y+\theta$ is an elliptic curve
$E^{\theta}$ with Weierstrass model \cite[Prop. 2.5]{nr}
\begin{equation}\label{Weierstrass}
E^{\theta}\colon\quad y^2+g\theta^{-1}y=x^3+dx^2+e.
\end{equation}
Hence, the Jacobian $J$ of $C$ is $\kb$-isogenous to the product
$E^{\theta}\times E^{\theta'}\times E^{\theta''}$ of the three supersingular elliptic curves corresponding to the three nonzero roots
$\theta,\,\theta',\,\theta''$ of the separable polynomial $R(Y)$.

In this section we determine the possible values of the isogeny classes of these triples of elliptic quotients of supersingular plane quartics. In section \ref{sec4} we shall deduce from these results the rational functions that can occur as the zeta function of a supersingular plane quartic.

We shall work from now on with the Artin-Schreier models of $E^{\theta}$, $E^{\theta'}$ and $E^{\theta''}$:
\begin{equation}\label{etheta}
E^{\theta}\colon \quad y^2+y=a_{\theta}x^3+b_{\theta}x^2+c_{\theta},
\end{equation}
with  $a_{\theta}=(g^{-1}\theta)^2,\ b_{\theta}=a_{\theta}d, \
c_{\theta}=a_\theta e$, and similar notations for $E^{\theta'}$ and
$E^{\theta''}$. Note that $a_{\theta}+a_{\theta'}+a_{\theta''}=0$.
We analyze separately the three cases that arise according to the
structure of the Galois set
$\{a_{\theta},\,a_{\theta'},\,a_{\theta''}\}$. To this purpose, the
following definition will be useful.

\begin{defi}
Let $C$ be a quartic given by an equation (\ref{quartic}). According to the number $3,\,1$ or $0$ of roots in $k^*$ of the polynomial $R(Y)=Y^4+fY^2+gY$, we say that $C$ is respectively split, of quadratic type or of cubic type.
\end{defi}

If $C$ is split then the Jacobian of $C$ is $k$-isogenous to the product of three elliptic curves, but the converse is not true (cf. section \ref{sec4}).

Throughout this section we shall make constant use of Theorem \ref{allcurves} and Proposition \ref{crit} without further mention.
If $q$ is a square we shall denote by $\enc$ ($\op{nc}$ for ``non-cube") any supersingular elliptic curve $k$-isomorphic to $E_a$ for some non-cube $a\in k^*$.
Thus, if $q$ is a square the five $k$-isogeny classes of supersingular elliptic curves are represented by $E_1$, $E_1'$, $\enc$, $\enc'$ and $E_0$.

\subsection{Quartics of cubic type}
If $C$ is a quartic of cubic type the elliptic curves $E^{\theta}$, $E^{\theta'}$, $E^{\theta''}$
are $k_3$-isogenous because they are Galois conjugate.

\begin{defi}
Let $E$ be a supersingular elliptic curve defined over $k_3$. We say that $E$ {\it is attained} if there exists a plane quartic of cubic type over $k$ whose Jacobian is $k_3$-isogenous to $E\times E\times E$, or equivalently, such that $E^{\theta}$ is $k_3$-isogenous to $E$.
\end{defi}

Our aim in this subsection is to prove the following result

\begin{teor}\label{teorsimple}
All supersingular elliptic curves $E$ defined over $k_3$ are
attained by a plane quartic of cubic type over $k$  with only two
exceptions: $E_1'$ is not attained over $\ff$ and $H$ is not
attained over $\ft$.
\end{teor}

The isogeny class of $E^{\theta}$ is computed in Proposition
\ref{crit} in terms of $a_{\theta}$, $b_{\theta}$, $c_{\theta}$. If
we let $C$ span all quartics of cubic type the parameters
$(a_{\theta},\,b_{\theta},\,c_{\theta})$ span all triples
$(a,\,ad,\,ae)$, where $a$ is an arbitrary element in $k_3^*$ with
$\tr_{k_3/k}(a)=0$, and $d,\,e$ belong to $k$. We have
$\tr_{k_3/k}(c_{\theta})=e\tr_{k_3/k}(a)=0$; thus, the coefficient
$e$ of the quartic does not play any role and we can assume that $e=0$.

From now on, the discussion of Theorem \ref{teorsimple} takes a
different form for $q$ square and $q$ non-square.

\subsubsection*{Case $q$ square}
Let us check first that $a$ can be a cube and a non-cube.
\begin{lem}\label{normcube}
An element $a\in k_3^*$ is a cube in $k_3^*$ if and only if  $N_{k_3/k}(a)$ is a cube in $k^*$.
\end{lem}

\begin{proof}
Clearly $a=u^3$ implies that $N_{k_3/k}(a)=N_{k_3/k}(u)^3$  is a cube in $k^*$. Conversely, suppose that $\rho:=N_{k_3/k}(a)$ is a cube in $k^*$; then it is necessarily a ninth power in $k_3^*$ because every element in $k^*$ is a cube in $k_3^*$. Since $q\equiv1\md3$ we have $a'=a^q=au^3$, and similarly $u'=u^q=uw^3$, for some $u,w\in k_3^*$. Hence
$$\rho=aa'a''=a(au^3)(a'(u')^3)=a^3u^9w^9
$$
and $a$ is a cube in $k_3^*$ because $a^3$ is a ninth power.
\end{proof}

\begin{cor}
There are cubes and non-cubes $a\in k_3^*$ with  $\tr_{k_3/k}(a)=0$.
\end{cor}

\begin{proof}
The map $k^*\lra k$ given by $x\mapsto x^2+x^{-1}$ cannot be onto. For any $j\in k$ which is not in the image of this map the polynomial $x^3+jx+1$ is irreducible and the roots of this polynomial are cubes in $k_3^*$. On the other hand, the roots of the irreducible polynomials of the type $x^3+j$ are non-cubes in $k_3^*$.
\end{proof}

\begin{prop}
The elliptic curves $E_0,\, E_1,\, \enc$ are attained.
\end{prop}

\begin{proof}
Take any $a\in k_3^*$ with $\tr_{k_3/k}(a)=0$. For $d=0$ we get $E^{\theta}\sim_{k_3}E_a$ and the two isogeny classes $E_1$, $\enc$ are attained, according to $a$ being a cube or not.

Take now $a=u^3$ for some $u\in k_3^*$. If $\sigma_1$, $\sigma_2$, $\sigma_3$ are the elementary symmetric functions of the conjugate elements $u$, $u'$, $u''$ we have
\begin{equation}\label{sigma}
0=u^3+(u')^3+(u'')^3=\sigma_1^3+\sigma_1\sigma_2+\sigma_3.
\end{equation}
Hence, $\tr_{k_3/k}(u)=\sigma_1\ne 0$. By duality relative to the extension $k/\ff$ (cf. Lemma \ref{onto}), there exists $d\in k$ such that $0\ne\tr_{k/\ff}(d\sigma_1)=\tr_{k_3/\ff}(du)$. Hence,
 $b_{\theta}u^{-2}=adu^{-2}=du\not\in\as^2(k_3)$ and $E^{\theta}\sim_{k_3}E_0$.
\end{proof}

In order to show that $E_1'$ (respectively $\enc'$) is attained we
need to find $a\in k_3$, $d\in k$ such that  $\tr_{k_3/k}(a)=0$, $a$
is a cube (respectively a non-cube) in $k_3$, and there is some $v\in
k_3$ satisfying $d=v+v^4a$ and $\tr_{k_3}(v^3a)=1$. This will be shown
in the following lemma:

\begin{lem}
If $q>4$ there exist $a,\,v\in k_3^*$ such that $\tr_{k_3/k}(a)=0$ and
\begin{enumerate}[a)]
\item $a$ is a cube (resp. $a$ is a non-cube) in $k_3^*$,
\item $v+v^4a$ belongs to $k$ and $\tr_{k_3}(v^3a)=1$.
\end{enumerate}
If $q=4$ this is still true if $a$ is a non-cube and false if $a$ is a cube.
\end{lem}

\begin{proof}
Take any $a\in k_3^*$ with $\tr_{k_3/k}(a)=0$. Let $r,\,s,\,t\in
k_3$ be determined by $a=r^8,\,a'=s^8,\,a''=t^8$. Note that
$r,\,s,\,t$ are conjugate over $k$. Take any $\la\in k$ and let
$x,\,y,\,z$ be the solution of $x+x^4a=y+y^4a'=z+z^4a''$ given in
(\ref{param1}) of Lemma \ref{parameterization}. Clearly, $x,\,y,\,z$
are conjugate  over $k$; in particular, the common value
$d:=x+x^4a=y+y^4a'=z+z^4a''$ belongs to $k$. We get all possible $v\in k_3$ such that $v+v^4a\in k$ by taking $v=x$ for all values of $\la$ in $k$. Moreover, by the
further computation of Lemma \ref{parameterization}
$$
\tr_{k_3}(x^3a)=\tr_k(x^3a+y^3a'+z^3a'')=\tr_k(A\la^9+B\la^3),
$$with $A=(rst)^{-1},\ B=(rst)^{-1}(r^7s+s^7t+t^7r)^{1/4}$.
To complete the proof we need only to show that the curve $\dd_{AB}$ given by the equation $y^2+y=Ax^9+Bx^3$ has less than $2q+1$ points. For $q>64$ this is guaranteed by Lemma \ref{dabc}, regardless of the fact that $a$ is a cube or a non-cube in $k_3^*$.

Let us discuss now the cases $q\le64$. By Lemma \ref{dabc} we still have $|\dd_{AB}(k)|<2q+1$ as long as $w<n$. To deal with the case $a$ non-cube assume that $a$ is the root of an irreducible polynomial $x^3+j$ for some $j\in k$. We have then $A=j^{-1/8}$, $B=0$ and  $w<n$ in all cases. In fact, $0=B\ne A^2=j^{-1/4}$ for $q=16$ and, since $j$ is not a cube in $k^*$, we have $A\not\in\fe$ for $q=64$ and $A+B=A\not\in\ft$ for $q=4$. To deal with the case $a$ cube assume that $a$ is the root of an irreducible polynomial $x^3+jx+1$ for some $j\in k$. We have then $A=1$, $B^8=j+j^{1/4}$ and $w<n$ if $q>4$; in fact, for $q=64$ we have $B\ne 0$ (otherwise $j\in\ff$ and $x^3+jx+1$ would split in $k^*$), and for $q=16$ we can take $j=1$ and $0=B\ne A^2$. Finally, for $q=4$ we have $B=0$, and the curve $\dd_{AB}$ has $9$ points in $k$; thus, $\tr_{k_3}(v^3a)=0$ for all $v\in k_3$ such that $v+v^4a$ belongs to $k$.
\end{proof}

\subsubsection*{Case $q$ non-square}
\begin{lem}
There exists $j\in k$ such that the polynomial $x^3+x^2+jx+(j+1)$ is irreducible. If $u\in k_3$ is a root of this polynomial then $\tr_{k_3/k}(u^3)=0$.
\end{lem}

\begin{proof}
The second statement is consequence of $\tr_{k_3/k}(u)=\tr_{k_3/k}(u^2)=1$. To show the existence of $j$, take any $a\in k_3^*$ with $\tr_{k_3/k}(a)=0$; since $q$ is non-square $a$ is a cube, say $a=u^3$, and $\tr_{k_3/k}(u)\ne0$ by (\ref{sigma}). Multiplying $u$ by $\tr_{k_3/k}(u)^{-3}$ we get an element $u$ with $\tr_{k_3/k}(u)=1$ and $\tr_{k_3/k}(u^3)=0$, and its minimal polynomial over $k$ is $x^3+x^2+jx+(j+1)$ for some $j\in k$.
\end{proof}

Taking $d=0$ we attain the curve $E_1$. In order to attain $H$ and $H'$ we need only to prove the following lemma:

\begin{lem}\label{nonsq}
Let $j\in k$ be such that the polynomial $x^3+x^2+jx+(j+1)$ is irreducible and let $u\in k_3$ be a root of this polynomial.

If $q>2$ for any $\eps\in\ft$ there exists $v\in k_3$ such that $(1+v+v^4)u^{-1}$ belongs to $k$ and $\tr_{k_3}(v^3+v)=\eps$.

If $q=2$ this is true for $\eps=1$ and false for $\eps=0$.
\end{lem}

\begin{proof}
Let $r,\,s,\,t\in k_3$ be determined by $u^{-1}=r^{16},\,(u')^{-1}=s^{16},\,(u'')^{-1}=t^{16}$. Note that $r,\,s,\,t$ are conjugate over $k$. Take any $\la\in k$ and let $x,\,y,\,z$ be the solution of $u^{-1}(1+x+x^4)=(u')^{-1}(1+y+y^4)=(u'')^{-1}(1+z+z^4)$ given in (\ref{param2}) of Lemma \ref{parameterization}. Clearly, $x,\,y,\,z$ are conjugate over $k$; in particular, the common value $d:=u^{-1}(1+x+x^4)=(u')^{-1}(1+y+y^4)=(u'')^{-1}(1+z+z^4)$ belongs to $k$. We get all possible $v\in k_3$ such that $(1+v+v^4)u^{-1}\in k$ by taking $v=x$ for all values of $\la$ in $k$. Moreover, by the further computation of Lemma \ref{parameterization}
$$
\tr_{k_3}(x^3+x)=\tr_k(x^3+x+y^3+y+z^3+z)=\tr_k(A\la^9+B\la^3+C\la+D),
$$with $A=R^{-21}$, $B=R^{-3}T^{1/4}$, $C=ST^{1/8}+R^3S^{16}$, $D=R^{24}S^{72}$, where
$$
R=rst, \ S=r^{-1}+s^{-1}+t^{-1}, \ T=r^{-6}s^{-42}+s^{-6}t^{-42}+t^6r^{-42}.
$$
To complete the proof we need only to show that the curve $\dd_{ABCD}$ given by the equation $y^2+y=Ax^9+Bx^3+Cx+D$ has more than $1$ point and less than $2q+1$ points.
 This is equivalent to $2\le|\dd_{ABC}(k)|\le2q$ for the curve $y^2+y=Ax^9+Bx^3+Cx$, since
$$
|\dd_{ABCD}(k)|=|\dd_{ABC}(k)|,\quad\mbox{or}\quad |\dd_{ABCD}(k)|+|\dd_{ABC}(k)|=2q+2,
$$according to $D$ belonging to $\as(k)$ or not. For $q>8$ the desired inequality is guaranteed by Lemma \ref{dabc}.

The reader can check that $T^8=j(j^2+j+1)(j^3+j+1)(j^6+j^5+j^4+j+1)$. For $q=8$ the polynomial $x^3+x^2+jx+j+1$ is irreducible only if $j^3+j^2+1=0$; hence, $B=R^{-3}T^{-1/4}\ne0$ and $2\le|\dd_{ABC}(k)|\le 2q$ by Lemma \ref{dabc}. For $q=2$ we have necessarily $j=0$, $A=C=D=1$, $B=0$ and the curve $\dd_{ABCD}$ has $1$ point in $k$ (the point at infinity); thus, $\tr_{k_3}(v^3+v)=1$ for all values of $v\in k_3$ such that $(1+v+v^4)u^{-1}$ belongs to $k$.
\end{proof}

\subsection{Quartics of quadratic type}\label{hs}
If $C$ is a quartic of quadratic type the elliptic curve $E^{\theta}$ is defined over $k$ whereas $E^{\theta'}$, $E^{\theta''}$ are defined over $k_2$ and Galois conjugate, so that they are
 $k_2$-isogenous.

\begin{defi}
 Let $(E,F)$ be a pair of supersingular elliptic curves defined respectively over $k$ and $k_2$. We say that the pair $(E,F)$ {\it is attained} if there exists a plane quartic of quadratic type over $k$ such that $E^{\theta}$ is $k$-isogenous to $E$ and the Jacobian is $k_2$-isogenous to $E\times F\times F$.
\end{defi}

Our aim in this subsection is to prove the following result

\begin{teor}\label{teorhalfsplit}
All pairs $(E,F)$ of supersingular elliptic curves, $E$ defined over $k$ and $F$ defined over $k_2$, are attained by a plane quartic of quadratic type over $k$, with the following exceptions:
\begin{itemize}
\item
The pairs $(E_1,E_1')$, $(E_1',E_1)$ are not attained over $\mathbb{F}_{16}$.
\item
The pairs $(H,E_1)$, $(H',E_1')$ are not attained over $\fe$.
\item
All pairs $(E,F)$ with $E\sim_k E_0,\,E_1,\,E_1'$ and $F\sim_{k_2} E_0,\,E_1,\,E_1'$, and the pairs $(E_1,\enc')$, $(E_1',\enc)$, $(\enc,E_1')$, $(\enc',E_1)$  are not attained over $\ff$.
\item Over $\ft$ the attained pairs are:  $(E_1,\enc)$, $(E_1,\enc')$, $(H,\enc')$, $(H',\enc)$.
\end{itemize}
\end{teor}

If we let $C$ span all quartics of quadratic type the parameters
$(a_{\theta},\,b_{\theta},\,c_{\theta})$,
$(a_{\theta'},\,b_{\theta'},\,c_{\theta'})$ span all triples
$(a,\,ad,\,ae)$, $(a',\,a'd,\,a'e)$,  where $a'$ is an arbitrary
element in $k_2\setminus k$, $a=\tr_{k_2/k}(a')$, and $d,\,e$ belong
to $k$.

\begin{lem}\label{symmetry}
If the pair $(E,F)$ is attained then the pair $(E',F')$ is attained.
\end{lem}

\begin{proof}
Suppose that $(E,F)$ is attained by a quartic with data $(a',d,e)$.
Since $a\ne 0$, by Lemma \ref{onto} we can always take $e_0\in k$ such that $\tr_k(ae_0)=1$. Since  $\tr_{k_2/k}(a'e_0)=ae_0$, we have in particular $\tr_{k_2}(a'e_0)=\tr_k(ae_0)=1$. Thus, the pair $(E',F')$ is attained by a quartic with data $(a',d,e+e_0)$.
\end{proof}

From now on, the discussion of Theorem \ref{teorhalfsplit} takes a
different form for $q$ square and $q$ non-square.

\subsubsection*{Case $q$ square}
\begin{lem}\label{cubenoncube}\mbox{\null}
\begin{enumerate}[a)]
\item If $q>4$, among the $k_2/k$-traces of cubes in $k_2\setminus k$
there are cubes and non-cubes of $k^*$. If $q=4$ the traces of cubes
in $k_2\setminus k$ are always non-cubes in $k^*$. \label{casea}
\item Among the
$k_2/k$-traces of non-cubes in $k_2\setminus k$ there are cubes and
non-cubes of $k^*$.
\end{enumerate}
\end{lem}

\begin{proof}
Take $z\in k\setminus \as(k)$ and $x\in k_2\setminus k$ such that $x^2+x=z+1$. We have $\tr_{k_2/k}(x^3)=\tr_{k_2/k}(x^2+(z+1)x)=z$.

The map $E_1\setminus\{\infty\}\lra \{x\in k\tq x^3\in\as(k)\}$ determined by $(x,y)\mapsto x$
is onto and $2\colon 1$; hence, the target set has cardinality $\left(|E_1(k)|-1\right)/2$.

Thus, there are $(q/2)+(-1)^{n/2}\sqrt{q}$ elements in $k$ whose cube is not in $\as(k)$ and there are $(q+(-1)^{n/2}2\sqrt{q})/6$ cubes in $k\setminus\as(k)$. This means that for $q>4$ there are elements $z$ as above being cubes and non-cubes.

For $q=4$, assume that $x^2+wx+z=0$ is the minimal polynomial of $x\in k_2\setminus k$. Then, $\tr_{k_2/k}(x^3)=1+wz$, and this trace is always a non-cube because $zw\ne0$. This proves the first item.

There are $2(q^2-1)/3-2(q-1)/3=2q(q-1)/3$ non-cubes in $k_2\setminus
k$. Their $k_2/k$-traces take at least $2(q-1)/3$ different values,
so that they cannot all be cubes. They cannot all be non-cubes
either. In fact, if they were all non-cubes, their traces would take
exactly $2(q-1)/3$ different values, and the set of non-cubes would
be closed under addition by elements of $k$. Therefore, the traces
of all cubes in $k_2\setminus k$ would all be cubes in $k$, in
contradiction with \ref{casea}).
\end{proof}

Taking quartics of quadratic type with $d=e=0$ we attain the pairs
$(E_1,\enc)$, $(\enc,E_1)$, $(\enc,\enc)$, and the pair $(E_1,E_1)$
if $q>4$, by considering the four possibilities for $a'$
cube/non-cube with $a=\tr_{k_2/k}(a')$ cube/non-cube.

\begin{prop}\label{bu}
If $q>4$ and $E\sim_kE_0$ or $F\sim_{k_2}E_0$, the pair $(E,F)$ is attained.

If $q=4$ the pairs $(E_0,\enc)$ and $(\enc,E_0)$ are attained.
\end{prop}

\begin{proof}
If $a'$ is a non-cube and $a=u^3$ for some $u\in k^*$, we can take $d\in k$ such that $\tr_{k/\ff}(du)\ne0$. We get $b_{\theta}u^{-2}=du\not\in\as^2(k)$  and we attain one of the pairs $(E_0,\enc)$ or $(E_0,\enc')$ (hence both are attained by Lemma \ref{symmetry}). Similarly, if $a'=(u')^3$ and $a$ is a non-cube, we attain one of the pairs
$(\enc,E_0)$ or $(\enc',E_0)$, just by taking $d\in k$ with $\tr_{k_2/\ff}(du')=\tr_{k/\ff}(d\tr_{k_2/k}(u'))\ne0$.

From now on we suppose $q>4$, $a'=(u')^3$ and $a=u^3$. We claim that
$u$ and $\la:=\tr_{k_2/k}(u')$ are $\ff$-linearly independent. In
fact, if $(u')^2=\la u'+\mu$ is the minimal equation of $u'$ over
$k$, we see that $$a=\tr_{k_2/k}((u')^3)=\tr_{k_2/k}(\la (u')^2+\mu
u')=\la^3+\mu\la.$$ Hence, $u=\omega\la$ for some $\omega\in\ff^*$
would imply $\la^3+\mu\la=a=u^3=\la^3$, which is impossible, since
$\la\mu\ne0$. By Lemma \ref{onto}, we can find $d\in k$ such that
$w:=\tr_{k/\ff}(du)$ and $w':=\tr_{k_2/\ff}(du')=\tr_{k/\ff}(d\la)$
take prescribed values. Therefore, the pair $(E_0,E_0)$ is attained
($ww'\ne0$), one of the pairs $(E_0,E_1)$ or $(E_0,E_1')$ is
attained ($w\ne0,\,w'=0$) and one of the pairs $(E_1,E_0)$ or
$(E_1',E_0)$ is attained ($w=0,\,w'\ne0$).
\end{proof}

In order to finish the proof of Theorem \ref{teorhalfsplit} in the case $q$ square, we have to check that the eight remaining pairs:
$$\ars{1.4}
\begin{array}{ll}
(E_1,E_1'),\ (E_1', E_1);&\quad (E_1,\enc'),\ (E_1', \enc);\\
(\enc,E_1'),\ (\enc', E_1);&\quad (\enc,\enc'),\ (\enc', \enc)
\end{array}
$$
are attained for $q>16$, the last six pairs for $q=16$, and only the last two pairs for $q=4$. By Lemma \ref{symmetry}, it is sufficient to prove the following result:

\begin{lem}\label{changesign}
Let $a'\in k_2\setminus k$ and let $a=\tr_{k_2/k}(a')$.

If $q>16$ there exist $d,\,v\in k$, $v'\in k_2$ such that $v+v^4a=d=v'+(v')^4a'$ and $\tr_k(v^3a)+\tr_{k_2}((v')^3a')=1$.

If $q=16$ this condition fails if and only if $a'$, $a$ are both cubes. If $q=4$ this condition holds if and  only if $a'$, $a$ are both non-cubes.
\end{lem}

\begin{proof}
Let $r\in k$, $s,\,t\in k_2$ be determined by
$a=r^8,\,a'=s^8,\,a''=t^8$. Note that $s,\,t$ are conjugate over $k$
and $r=s+t$. Take any $\la\in k$ and let $x,\,y,\,z$ be the solution
of $x+x^4a=y+y^4a'=z+z^4a''$ given in (\ref{param1}) of Lemma
\ref{parameterization}. Clearly, $x$ belongs to $k$ and $y,\,z$ are
conjugate over $k$; in particular, the common value
$d:=x+x^4a=y+y^4a'=z+z^4a''$ belongs to $k$. We get all pairs $v\in k$, $v'\in k_2$ such that $v+v^4a=v'+(v')^4a'$ by taking $v=x$, $v'=y$ for all values of $\la$ in $k$. By the further
computation of Lemma \ref{parameterization}
$$
\tr_k(x^3a)+\tr_{k_2}(y^3a)=\tr_k(x^3a+y^3a'+z^3a'')=\tr_k(A\la^9+B\la^3),
$$with $A=(rst)^{-1},\ B=(rst)^{-1}(r^7s+s^7t+t^7r)^{1/4}$.
To complete the proof we need only to show that the curve $\dd_{AB}$
given by the equation $y^2+y=Ax^9+Bx^3$ has less than $2q+1$ points.
For $q>64$ this is guaranteed by Lemma \ref{dabc}.

By Lemma \ref{dabc}, for $q\le64$ the curve $\dd_{AB}$ has less than $2q+1$ points if and only if $w<n$. We check now when the conditions on $A,\,B$ ensuring that $w<n$ are satisfied. The minimal equation of $s,\,t$ over $k$ is $x^2+rx+st$; thus, $u:=r^{-2}st$ does not belong to $\as(k)$ (cf. Lemma \ref{qeq}). Let us express $A,\,B$ in terms of $u$: $A=(r^3u)^{-1}$, $B^4=(ru)^{-4}(u^4+u^3+u^2+1)=(ru)^{-4}(u+1)(u^3+u+1)$.

For $q=64$ we have $B\ne0$; in fact, $B=0$ leads to $u\in \fe\subset\as(k)$.
For $q=16$, the condition $w=n$ is equivalent to $B=A^2$; raising to
the fourth power, this is equivalent to $u^4(u^4+u^3+u^2+1)=r^{-5}$.
The eight elements of $k=\mathbb{F}_{16}$ that are not in $\as(k)$
are the four roots of $x^4+x^3+1=0$ and the four roots of
$x^4+x^3+x^2+x+1=0$, both polynomials irreducible over $\ft$.  Now,
if $u^4+u^3=1$ we have $B\ne A^2$ because $B=A^2$ implies
$u^6=r^{-5}$, which leads to a contradiction: $1=u^{18}=u^3$.
Therefore, if $u^4+u^3+u^2+u=1$, the condition $B=A^2$ is equivalent
to $u^5=1=r^{-5}$, which is equivalent to $u,\,r$ being both cubes in $k^*$. It is easy to check that this is equivalent to $a,\,a'$ being both cubes.
For $q=4$ the condition $A+B\in\ft$ translates into $r=1$ or $r=u$.
Hence, $A+B\in\ft$ if $a$ is a cube ($r=1$) or $a'$ is a cube ($st=1,\,u=r^{-2}=r$), and $A+B\not\in\ft$ if $a,\,a'$ are non-cubes ($r\ne1,\,st\ne 1,\,u=rst\ne r$).
\end{proof}

\subsubsection*{Case $q$ non-square}
For $q=2$ it is easy to check that the only attained pairs are:
$$
(E_1,\enc),\ (E_1,\enc'),\ (H,\enc'),\ (H',\enc).$$
From now on we
assume $q>2$. We can assume also that $a=1$; in fact, a quartic with
relevant data $(a',d,e)$ has the same elliptic quotients
$E^{\theta}$, $E^{\theta'}$ than a quartic with data $(a'/a,ud,ae)$,
where $a=\tr_{k_2/k}(a')$ and $a=u^3$.

\begin{prop}
If $E\sim_kE_1$ or $F\sim_{k_2}E_0$, the pair $(E,F)$ is attained.
\end{prop}

\begin{proof}
We attain $(E_1,E_1)$ and $(E_1,\enc)$  by taking $d=e=0$ and $a'$
cube or non-cube respectively. By Lemma \ref{symmetry} the pairs $(E_1,E_1')$, $(E_1,\enc')$ are attained too.

Let $a'=(u')^3$ be a cube with $\tr_{k_2/k}(a')=1$. Arguing as in Proposition \ref{bu}, we see that $1\ne\la:=\tr_{k_2/k}(u')$ and by Lemma \ref{onto}, we can find $d\in k$ such that $w:=\tr_k(d)$ and $w':=\tr_{k_2}(du')=\tr_k(d\la)$ take prescribed values. Note that $w'=1$ ensures in particular that $\tr_{k_2/\ff}(du')\ne0$.
Therefore, the pair $(E_1,E_0)$ is attained ($w=0,\,w'=1$), and one of the pairs $(H,E_0)$ or $(H',E_0)$ (hence both by Lemma \ref{symmetry}) is attained ($w=w'=1$) .
\end{proof}

In order to finish the proof of Theorem \ref{teorhalfsplit} in the case $q$ non-square, we have to check that the eight remaining pairs:
$$\ars{1.4}
\begin{array}{ll}
(H,E_1),\ (H', E_1');&\quad (H,E_1'),\ (H', E_1);\\
(H,\enc),\ (H', \enc');&\quad (H,\enc'),\ (H', \enc)
\end{array}
$$
are attained, the first two pairs only for $q>8$. By Lemma \ref{symmetry}, it is sufficient to prove the following result:

\begin{lem}
Let $a'\in k_2$ with $\tr_{k_2/k}(a')=1$ and let $\eps\in\ft$.

If $q>8$ there exist $v\in k$, $v'\in k_2$ such that $1+v+v^4=v'+(v')^4a'$ and $\tr_k(v^3+v)+\tr_{k_2}((v')^3a')=\eps$.

If $q=8$ this is still true if $a'$ is a non-cube or $\eps=1$, and false otherwise.
\end{lem}

\begin{proof}
Let $a''$ be the conjugate of $a'$ over $k$ and write $a'=s^8$,
$a''=t^8$. Note that $s,\,t$ are conjugate over $k$ and $s+t=1$. We
can find a parametrization of the spatial curve
$1+x+x^4=y+y^4a'=z+z^4a''$ as we did in Lemma
\ref{parameterization}:
$$\ars{1.6}
\begin{array}{l}
x=1+(st)^{-1}+(st)^{-2}+(st)^{-1}\la+(1+(st)^{-1}+(st)^{-2})\la^4+(st)^{-2}\la^{16},\\
y=s^{-1}t^{-2}+s^{-4}t^2+(st)^{-1}\la+(s^{-1}t^{-2}+s^{-4}t^2)\la^4+t^{-2}s^{-4}\la^{16},\\
z=t^{-1}s^{-2}+t^{-4}s^2+(st)^{-1}\la+(t^{-1}s^{-2}+t^{-4}s^2)\la^4+s^{-2}t^{-4}\la^{16}.
\end{array}
$$
If we choose $\la\in k$ we get $x\in k$ and $y,\,z\in k_2$ conjugate over $k$. Let us denote
$u:=st$; note that $u\not\in \as(k)$ because the polynomial $x^2+x+st$ is irreducible over $k$. By a straightforward computation
\begin{multline*}
x^3+x+y^3a'+z^3a''=u^{-4}\la^{36}+u^{-4}\la^{32}+u^{-2}\la^{18}+u^{-2}\la^{16}+\\+(1+u^{-1}+u^{-2}+u^{-4})\la^{12}+u^{-1}\la^9
+(1+u^{-1}+u^{-2}+u^{-4})\la^8+u^{-4}\la^4+u^{-4}.
\end{multline*}
Therefore, $\tr_k(x^3+x+y^3a'+z^3a'')=\tr_k(A\la^9+B\la^3+C\la+D)$, with $A=u^{-1}$, $B=1+u^{-1/4}+u^{-1/2}+u^{-1}$, $C=1+u^{-1/8}+u^{-1/4}+u^{-1/2}+u^{-1}$, $D=u^{-4}$.
To complete the proof we need only to show that the curve $\dd_{ABCD}$ given by the equation $y^2+y=Ax^9+Bx^3+Cx+D$ has more than $1$ point and less than $2q+1$ points. Arguing as in Lemma \ref{nonsq}, this is equivalent to $2\le|\dd_{ABC}(k)|\le2q$ for the curve $y^2+y=Ax^9+Bx^3+Cx$. For $q>8$ this is guaranteed by Lemma \ref{dabc}.

For $q=8$ this is still true if $B\ne0$.  Now, $B^4=u^{-4}(u+1)(u^3+u+1)$, so that $B=0$ if and only if $u=1$, or equivalently $a'\in\ff\setminus\ft$; in this case $a',\,a''$ are the only cubes in $k_2$ with relative trace $1$ over $k$. Finally, in this bad case we have $B=0$, $A=C=D=1$, the curve  $\dd_{ABCD}$ has
only one rational point and $\tr_k(v^3+v)+\tr_{k_2}((v')^3a')=1$ for all possible values of $v,\,v'$.
\end{proof}

\subsection{Split quartics}
\begin{defi}
We say that an unordered triple $\{E,\,F,\,G\}$ of supersingular curves defined over $k$ {\it  is attained} if there exists a split plane quartic over $k$ whose Jacobian is $k$-isogenous to $E\times F\times G$. We think $\{E,\,F,\,G\}$ as a multiset, allowing repetitions.

The bitwists of a triple $\{E,\,F,\,G\}$ are the three triples obtained by twisting exactly two of the curves $E,\,F,\,G$.
\end{defi}

There is no split quartic over $\ft$; thus, we assume $q>2$ throughout this subsection. Our aim is to prove the following result

\begin{teor}\label{teorsplit}
All triples $\{E,\,F,\,G\}$ of supersingular elliptic curves over $k$ are attained by a split quartic over $k$, with the following exceptions:

\begin{itemize}
\item The triple $\{E_1,\,E_1,\,E_1'\}$ and its bitwists are not attained over $\mathbb{F}_{64}$.
\item The triple $\{E_1,\,E_1,\,\enc'\}$ and its bitwists as well as
 the triples $\{E,F,G\}$ with $E,\,F,\,G\in\{E_1,E_1',E_0\}$ are not attained over
$\mathbb{F}_{16}$.
\item The triple $\{H,\,H,\,H\}$ and its bitwists are not attained over $\fe$.
\item For $k=\ff$ the only triples attained are $\{E_0,\,\enc,\,\enc\}$, $\{E_1,\,\enc,\,\enc\}$ and their bitwists.
\end{itemize}
\end{teor}

Given a split quartic $C$, the relevant data to determine the
$k$-isogeny classes of the three elliptic quotients of $C$ are the
five elements  $a,\,a',\,a'',\,d,\,e\in k$, where $a=a_{\theta}$,
$a'=a_{\theta'}$, $a''=a_{\theta''}$. If $C$ spans all split
quartics the set $\{a,a',a''\}$ spans all triples in $k^*$
satisfying $a+a'+a''=0$.

\begin{lem}\label{change2signs}
If the triple $\{E,\,F,\,G\}$ is attained then all its bitwists are attained.
\end{lem}

\begin{proof}
Suppose the Jacobian $J$ of a quartic with parameters
$(a,a',a'',d,e)$ satisfies $J\sim_kE\times F\times G$. By duality,
there exists $e_0\in k$ such that $\tr_k(ae_0)=\tr_k(a'e_0)=1$ (and
$\tr_k(a''e_0)=0$). Hence, the Jacobian of the quartic with
parameters $(a,a',a'',d,e+e_0)$ is $k$-isogenous to $E'\times
F'\times G$.
\end{proof}

From now on, the discussion of Theorem \ref{teorsplit} takes a
different form for $q$ square and $q$ non-square.

\subsubsection*{Case $q$ square}
\begin{lem}\label{yyy}
Given any natural number $0\le i\le 3$ there are elements $x,y,z\in k^*$ such that exactly $i$ of them are cubes and $x+y+z=0$, except in the following cases:

If $x,y,z$ are cubes in $\mathbb{F}_{16}^*$ then $x+y+z\ne0$.

If $x,y,z\in\ff^*$ satisfy $x+y+z=0$ then necessarily two of the elements $x,y,z$ are non-cubes and the third is a cube.
\end{lem}

\begin{proof}
By Lemma \ref{surface} there are three different cubes in $k^*$ adding up to zero if and only if $q>16$. If $q>4$, there are two different cubes whose sum is a non-cube, because $(k^*)^3\cup\{0\}$ is not an additive subgroup of $k$ (the cardinality is not a divisor of $q$). The same argument shows that for all $q$ there are two different non-cubes whose sum is a cube. Finally, let $\op{NC}\subseteq k^*$ be the set of non-cubes of $k^*$; for any $x\in \op{NC}$ the set $x+\op{NC}$ contains
$2(q-1)/3$ different elements so that it must contain a non-cube if $q>4$.
\end{proof}

\begin{prop}\label{ezero}
If $q>16$ the triples $\{E_0,\,F,\,G\}$ are all attained.

If $q=16$ the triple $\{E_0,\,F,\,G\}$ is attained if and only if $F,\,G$ don't belong both to $\{E_1,E_1',E_0\}$.

If $q=4$ the triple $\{E_0,\,F,\,G\}$ is attained if and only if $F,\,G$ are both $k$-isogenous to $\enc$ or $\enc'$.
\end{prop}

\begin{proof}
If $q>16$ there is a solution in $k$ to $x^3+y^3+z^3=0$, $xyz\ne0$. We consider a quartic with parameters $a=x^3,\,a'=y^3,\,a''=z^3$. We claim that the elements $x,\,y,\,z$ are $\ff$-linearly independent. In fact, any two of them are not $\ff$-proportional, since $y=\la x$ with $\la\in\ff^*$, implies $y^3=x^3$, leading to $z=0$, against our assumption. Also, if we had $z=\la x+\mu y$ with $\la,\,\mu\in\ff^*$ we would have $x^3+y^3=z^3=x^3+y^3+\la\mu xy(\la x+\mu y)$, leading to $\la\mu xy(\la x+\mu y)=0$, which is impossible.

By duality, we can find $d\in k$ such that $\tr_{k/\ff}(dx)$, $\tr_{k/\ff}(dy)$, $\tr_{k/\ff}(dz)$ takes prescribed values. Hence we can choose $\tr_{k/\ff}(dx)$ to be nonzero in order to achieve $E^{\theta}\sim_kE_0$ and, having in mind Lemma \ref{change2signs}, we can obtain a quartic attaining any triple $\{E_0,\,F,\,G\}$ with $F,\,G\in\{E_0,\,E_1,\,E_1'\}$.

If $q>4$ we can attain any other triple containing $E_0$ by the same (even easier) argument: we start with any solution in $k$ of $a+a'+a''=0$
with $a=x^3$, $a''$ non-cube and $a'$ cube/non-cube at will. We choose a parameter $d \in k$ such that $dx\not\in \as^2(k)$ and $dy$ belonging or not to $\as^2(k)$  at will, if $a'=y^3$. Finally, with a right choice of the parameter $e$ we can perform two arbitrary twists on $E_0,\,F,\,G$ by Lemma \ref{change2signs}, which in practice means that we can perform arbitrary twists on $F$, $G$, since $E_0'=E_0$.

For $q=4$ we can do exactly the same, but restricted to the situation in which two exactly of $a$, $a'$, $a''$ are non-cubes, as Lemma \ref{yyy} shows.
\end{proof}

We can attain the following triples which do not contain $E_0$, just
by taking $d=e=0$, and $a+a'+a''=0$ being cubes or non-cubes at
will:
$$
\{E_1,\,E_1,\,E_1\},\quad\{E_1,\,E_1,\,\enc\},\quad\{E_1,\,\enc,\,\enc\},\quad\{\enc,\,\enc,\,\enc\},
$$the first triple only if $q>16$, the second and fourth only if $q>4$.
By Lemma \ref{change2signs} we can attain all their bitwists too. To finish the proof of Theorem \ref{teorsplit} we need only to attain the triples that are obtained by applying three twists (or equivalently one twist) to the four triples in the above list. This will be a consequence of the following lemma:

\begin{lem}
Suppose that $a,\,a',\,a''\in k^*$ satisfy $a+a'+a''=0$. There exist $x,\,y,\,z\in k$ such that
\begin{equation}\label{spatial0}
x+x^4a=y+y^4a'=z+z^4a''
\end{equation} and $\tr_k(x^3a+y^3a'+z^3a'')=1$, except in the following cases:
\begin{itemize}
\item $q=64$ and $a,\,a',\,a''$ are all cubes in $k^*$,
\item $q=16$ and among $a,\,a',\,a''$ there are two cubes and one non-cube in $k^*$,
\item $q=4$.
\end{itemize}
\end{lem}

\begin{proof}
Let $r,\,s,\,t\in k$ be determined by $a=r^8,\,a'=s^8,\,a''=t^8$. Take
any $\la\in k$ and let $x,\,y,\,z$ be the solution of (\ref{spatial0})
 given in (\ref{param1}) of Lemma \ref{parameterization}. We want to
show the existence of $\la\in k$ such that
$$
1=\tr_k(x^3a+y^3a'+z^3a'')=\tr_k(A\la^9+B\la^3),
$$with $A=(rst)^{-1},\ B=(rst)^{-1}(r^7s+s^7t+t^7r)^{1/4}$ (cf. Lemma \ref{parameterization}).
To complete the proof we need only to show that the curve $\dd_{AB}$ given by the equation $y^2+y=Ax^9+Bx^3$ has less than $2q+1$ points. For $q>64$ this is guaranteed by Lemma \ref{dabc}.

For $q\le64$ the condition $|\dd_{AB}(k)|<2q+1$ is equivalent to $w<n$. We saw in the proof of Lemma \ref{changesign} that we can express $A=(r^3u)^{-1}$, $B^4=(ru)^{-4}(u+1)(u^3+u+1)$ in terms of $u:=r^{-2}st$, which now belongs to $\as(k)$.

For $q=64$ the condition $w=n$ is equivalent to $A\in\fe$, $B=0$,
which is equivalent to $r,\,s,\,t$ being all cubes in $k^*$. In
fact, if  $A\in\fe$ and $B=0$ we have $u,\,r^3\in\fe$ so that
$r^{21}=1$ and $r$ is a cube in $k^*$; since $A$ and $B$ are
symmetric in $r,\,s,\,t$ the other two elements $s,\,t$ are cubes
too. Conversely, if $r,\,s,\,t$ are cubes, they are all of the form
$\la\mu$, with $\la\in\ff^*$, $\mu\in\fe^*$; we can deduce from
Lemma \ref{surface} that there are two possibilities up to
permutation of $r,\,s,\,t$: either $(r,s,t)=(1,\omega,\omega^2)\mu$ for some $\omega\in\ff\setminus\ft$,
in which case $u=1$, $A=\mu^{-3}$, or
$(r,s,t)=\la(\mu_1,\mu_2,\mu_3)$, with $\mu_1+\mu_2+\mu_3=0$ in
$\fe$, in which case $u=(\mu_2/\mu_1)^2+(\mu_2/\mu_1)\in \as(\fe)$
(hence $u^3+u+1=0$) and $A=(\mu_1\mu_2\mu_3)^{-1}$.

For $q=16$, the condition $w=n$ is equivalent to $u^4(u^4+u^3+u^2+1)=r^{-5}$. This is equivalent to $u^4+u+1=0$ or $u=r^{-5}$ according to $r$ being a cube or not. The reader can check that this happens if and only if among $a,\,a',\,a''$ there are two cubes and one non-cube.

For $q=4$ we have necessarily $A=1$, $B=0$ and $w=n$ by Lemma \ref{dabc}.
\end{proof}

\subsubsection*{Case $q$ non-square}
We recall that we are assuming $q>2$. An argument completely analogous to that of Proposition \ref{ezero} shows that
\begin{prop}
The triples $\{E_1,\,F,\,G\}$ are all attained.
\end{prop}

Thus, in order to finish the proof of Theorem \ref{teorsplit} we
have only to check that all triples that can be built using the
curves $H$, $H'$ are attained, with some exceptions for $q=8$.
Having in mind Lemma \ref{change2signs} this will be a consequence
of the following lemma:

\begin{lem}
Let $u,\,u',\,u''\in k^*$ such that $u^3+(u')^3+(u'')^3=0$, and let $\eps\in\ft$.
If $q>8$ there exist $x,\,y,\,z\in k$ such that
\begin{equation}\label{spatial}
(1+x+x^4)u^{-1}=(1+y+y^4)(u')^{-1}=(1+z+z^4)(u'')^{-1}\end{equation}
and $\tr_k(x^3+x+y^3+y+z^3+z)=\eps$.

If $q=8$ the same is true for $\eps=1$ and false for $\eps=0$.
\end{lem}

\begin{proof}
Let $r,\,s,\,t\in k^*$ be determined by
$u^{-1}=r^{16},\,(u')^{-1}=s^{16},\,(u'')^{-1}=t^{16}$. Take any
$\la\in k$ and let $x,\,y,\,z$ be the solution of (\ref{spatial})
given in (\ref{param2}) of Lemma \ref{parameterization}. By the
further computation of Lemma \ref{parameterization} $$
\tr_k(x^3+x+y^3+y+z^3+z)=\tr_k(A\la^9+B\la^3+C\la+D),
$$with $A=R^{-21}$, $B=R^{-3}T^{1/4}$, $C=ST^{1/8}+R^3S^{16}$, $D=R^{24}S^{72}$, where
$$
R=rst, \ S=r^{-1}+s^{-1}+t^{-1}, \ T=r^{-6}s^{-42}+s^{-6}t^{-42}+t^{-6}r^{-42}.
$$
To complete the proof we need only to show that the curve $\dd_{ABCD}$ given by the equation $y^2+y=Ax^9+Bx^3+Cx+D$ has more than $1$ point and less than $2q+1$ points. For $q>8$ this is guaranteed by Lemma \ref{dabc}.

Assume $q=8$. We have $A=1$, and $B=0$ because $T=r^{-6}+s^{-6}+t^{-6}=0$. Note that $r+s+t=0$ because $x=x^{-6}$ for all $x\in k^*$; using this  one checks easily that $C=D=1$. Hence, the curve $\dd_{ABCD}$ has only one rational point.
\end{proof}

\section{Jacobians in isogeny classes of supersingular abelian threefolds}\label{sec4}
We want to find all possible values of the Weil polynomials of the
Jacobians of supersingular plane quartics $C$ defined over $k$, and
compare the results with the possible values for supersingular
abelian threefolds over $k$. Note that we do not need to consider
hyperelliptic curves since  Oort proved that there are no
supersingular hyperelliptic curves of genus three in characteristic
$2$ \cite[\S 5]{oort}.

If the quartic is split the Jacobian is $k$-isogenous to
$E^{\theta}\times E^{\theta'}\times E^{\theta''}$. If the quartic is
not split over $k$, the following result shows that the Weil
polynomial is also determined by the isogeny classes of the elliptic
quotients $E^{\theta}$, $E^{\theta'}$, $E^{\theta''}$ over their
field of definition.

\begin{lem}\label{ck}\mbox{\null}

(1) If $C$ is of cubic type then $|C(k)|=q+1$ and the Weil polynomial of the Jacobian is $x^6+tx^3+q^3$, where $x^2+t+q^3$ is the Frobenius polynomial of the elliptic quotient $E^{\theta}$ over $k_3$.

(2) If $C$ is of quadratic type then $|C(k)|=|E^{\theta}(k)|$, where $E^{\theta}$ is the elliptic quotient of $C$ defined over $k$. The Weil polynomial of the Jacobian is $(x^2+sx+q)(x^4+tx^2+q^2)$, where $x^2+sx+q$, $x^2+tx+q^2$ are the respective Frobenius polynomial of $E^{\theta}/k$ and $E^{\theta'}/k_2$.
\end{lem}

\begin{proof}
Let $C$ be given by  (\ref{quartic}). The polynomial
$R(Y)=Y^4+fY^2+gY$ determines an $\ft$-linear map $R\colon k\lra k$.
If $R(Y)$ has only the root $Y=0$ in $k$, this linear map is an
isomorphism and for each value $X\in k$ there is a unique value
$Y\in k$ such that $(X,Y)$ belongs to $C(k)$. This shows that
$|C(k)|=q+1$, $|C(k_2)|=q^2+1$, and these conditions imply that the
Weil polynomial over $k$ of the Jacobian $J$ of $C$ is of the type
$x^6+tx^3+q^3$ for some integer $t$. One checks easily that the Weil polynomial of $J$
over $k_3$ is then $(x^2+tx+q^3)^3$. Since
$J\sim_{k_3}E^{\theta}\times E^{\theta}\times E^{\theta}$, the
polynomial $x^2+tx+q^3$ is necessarily the Frobenius polynomial of
$E^{\theta}$. This proves the first item.

If $R(Y)$ has only one nonzero root $\theta$ in $k$ the linear map $Q\colon k\lra k$ determined by $Q(y)=y^2+g\theta^{-1}y$ has the same image than $R$. In fact,
these subspaces have both codimension one in $k$ because $\k(R)=\{0,\theta\}$ and $\k(Q)=\{0,g\theta^{-1}\}$; on the other hand, $R(Y)=Q(Y^2+\theta Y)$ (since $\theta^2+g\theta^{-1}=f$), so that $R(k)\subseteq Q(k)$ and both subspaces coincide. Therefore, for each $x\in k$, both $C$ and $E^{\theta}$ have two or zero $k$-rational points
having $x$ as first coordinate, according to $x^3+dx^2+e$ belonging or not to $R(k)=Q(k)$.

Let $x^2+sx+q$ be the Frobenius polynomial of $E^{\theta}/k$ and let $J\sim_k E^{\theta}\times A$ for some abelian surface $A$. The condition $|C(k)|=|E^{\theta}(k)|$ implies that  the Weil polynomial of $J$ over $k$ is of the type $(x^2+sx+q)(x^4+tx^2+q^2)$. One checks easily that the Weil polynomial of $A$ over $k_2$ is $(x^2+tx+q^2)^2$. Since $J\sim_{k_2}E^{\theta}\times E^{\theta'}\times E^{\theta'}$, the polynomial $x^2+tx+q^2$ is necessarily the Frobenius polynomial of $E^{\theta'}$.
\end{proof}

Let us display the isogeny classes of simple supersingular abelian varieties of dimension $1,\,2,\,3$. The Frobenius polynomials of the supersingular elliptic curves over $k$ are given in Proposition \ref{allcurves}. Let $A_{(s,t)}$ denote the isogeny class of the abelian surfaces over $k$ with Weil polynomial $x^4+sx^3+tx^2+qsx+q^2$. The $k$-simple supersingular isogeny classes of abelian surfaces are \cite[Theorem 2.9]{mn02}:
$$
A_{(0,0)},\quad A_{(0,-q)},\quad A_{(\sqrt{q},q)},\quad A_{(-\sqrt{q},q)},\qquad\qquad\mbox{ if $q$ is a square, and}
$$
$$
A_{(0,-2q)},\quad A_{(0,q)},\quad A_{(0,-q)},\quad A_{(\sqrt{2q},q)},\quad A_{(-\sqrt{2q},q)},\qquad\mbox{ if $q$ is not a square.}
$$
Finally, it is easy to deduce from a
result of Zhu \cite[Prop. 3.1]{zhu} that there is no simple
supersingular abelian threefold if q is nonsquare, and for q square
there are only two isogeny classes corresponding to the polynomials
$ x^6\pm \sqrt{q^3}x^3+q^3$.

Just by checking this list, Lemma \ref{ck} tells us when the isogeny class of the Jacobian of a supersingular quartic of cubic or quadratic type is simple, an elliptic curve times a simple surface or the product of three elliptic curves, in terms of the elliptic quotients of the quartic. The results are displayed in Tables 5,6,7 and 8.

\begin{table}
\begin{center}
\renewcommand{\arraystretch}{1.4}
\begin{tabular}{|c|c|c|}
\hline $E^{\theta}/k_3$&Weil polynomial of $J$&$k$-isogeny class of $J$\\
\hline$E_0$&$x^6+q^3$\qquad &$E_0\times A_{(0,-q)}$\\
\hline$E_1$&$ x^6-2\eps\sqrt{q^3}x^3+q^3$&$E_1\times \enc\times\enc$\\
\hline$E_1'$&$x^6+2\eps\sqrt{q^3}x^3+q^3$&$E_1'\times\enc'\times\enc'$\\
\hline$\enc,\,\enc'$&$x^6\pm\sqrt{q^3}x^3+q^3$&$k$-simple\\
\hline\end{tabular}\vskip.2cm
\caption{Isogeny class of the Jacobian of a quartic of cubic type, $q$ square, $\eps=(-1)^{n/2}$}
\label{Table:JcubicSq}
\end{center}
\end{table}
\begin{table}
\begin{center}
\renewcommand{\arraystretch}{1.4}
\begin{tabular}{|c|c|c|}
\hline $E^{\theta}/k_3$&Weil polynomial of $J$&$k$-isogeny class of $J$\\
\hline$E_1$&$x^6+q^3$ &$E_1\times A_{(0,-q)}$\\
\hline$H$&$x^6-\eps\sqrt{2q^3}x^3+q^3$&$H\times A_{(-\eps\sqrt{2q},q)}$\\
\hline$H'$&$x^6+\eps\sqrt{2q^3}x^3+q^3$&$H'\times A_{(\eps\sqrt{2q},q)}$\\
\hline\end{tabular}\vskip.2cm
\caption{Isogeny class of the Jacobian of a quartic of cubic type, $q$ non-square, $\eps=(-1)^{(n^2-1)/8}$}
\label{Table:JcubicNsq}
\end{center}
\end{table}
\begin{table}
\begin{center}
\renewcommand{\arraystretch}{1.4}
\begin{tabular}{|c|c|c|}
\hline $E^{\theta'}/k_2$&Weil polynomial of $J/E^{\theta}$&$k$-isogeny class of $J$\\
\hline$E_0$&$x^4+q^2$ &$E^{\theta}\times A_{(0,0)}$\\
\hline$E_1$&$x^4-2qx^2+q^2$&$E^{\theta}\times E_1\times E_1'$\\
\hline$E_1'$&$x^4+2qx^2+q^2$&$E^{\theta}\times E_0\times E_0$\\
\hline$\enc$&$x^4+qx^2+q^2$&$E^{\theta}\times\enc\times\enc' $\\
\hline$\enc'$&$x^4-qx^2+q^2$&$E^{\theta}\times A_{(0,-q)} $\\
\hline\end{tabular}\vskip.2cm
\caption{Isogeny class of the Jacobian of a quartic of quadratic type, $q$ square, $E^{\theta}$ is the elliptic quotient defined over $k$}
\label{Table:JquadraticSq}
\end{center}
\end{table}
\begin{table}
\begin{center}
\renewcommand{\arraystretch}{1.4}
\begin{tabular}{|c|c|c|}
\hline $E^{\theta'}/k_2$&Weil polynomial of $J/E^{\theta}$&$k$-isogeny class of $J$\\
\hline$E_0$&$x^4+q^2$&$E^{\theta}\times H\times H'$\\
\hline$E_1$&$x^4+2qx^2+q^2$&$E^{\theta}\times E_1\times E_1$\\
\hline$E_1'$&$x^4-2qx^2+q^2$&$E^{\theta}\times A_{(0,-2q)}$\\
\hline$\enc$&$x^4-qx^2+q^2$&$E^{\theta}\times A_{(0,-q)} $\\
\hline$\enc'$&$x^4+qx^2+q^2$&$E^{\theta}\times A_{(0,q)} $\\
\hline\end{tabular}\vskip.2cm
\caption{Isogeny class of the Jacobian of a quartic of quadratic type, $q$ non-square, $E^{\theta}$ is the elliptic quotient defined over $k$}
\label{Table:JquadraticNsq}
\end{center}
\end{table}

Theorems \ref{teorsimple}, \ref{teorhalfsplit} and \ref{teorsplit}
determine the existence of supersingular quartics with elliptic quotients in prescribed isogeny classes. These results, thanks to Tables 5,6,7,8, can now
be reinterpreted as the determination of all polynomials that occur as the Weil polynomial of the Jacobian of a supersingular quartic, giving thus a complete picture of which isogeny classes of supersingular abelian
threefolds contain a Jacobian.

\begin{teor} \label{thsimple}
All isogeny classes of simple supersingular abelian threefolds over
$k$ contain Jacobians.
\end{teor}

\begin{teor} \label{thhalf}
Let $f(x)$ be the Frobenius polynomial of a supersingular elliptic
curve over $k$ and let $g(x)$ be the Frobenius polynomial of a
supersingular $k$-simple abelian surface over $k$. Then, the isogeny
class of abelian threefolds with Weil polynomial $f(x)g(x)$ contains
a Jacobian, except in the following cases (where $\eps=\pm1$ takes
both values $1$ and $-1$):
\begin{itemize}
\item
$q$ is a square and $g(x)=x^4+\eps\sqrt{q}x^3+qx^2+\eps q\sqrt{q}x+q^2$,
\item
$q=4$ and $f(x)g(x)$ is one of the following polynomials
$$
(x^2+4)(x^4+16),\quad(x^2\pm 4x+4)(x^4+16),\quad (x^2+4x+4)(x^4-4x^2+16),
$$
\item
$q$ is a non-square and  $f(x)g(x)$ is one of the following polynomials
$$\ars{1.6}
\begin{array}{l}
(x^2+q)(x^4+\eps\sqrt{2q}x^3+qx^2+\eps q\sqrt{2q}x+q^2),\\
(x^2+\eps\sqrt{2q}x+q)(x^4+\eps\sqrt{2q}x^3+qx^2+\eps q\sqrt{2q}x+q^2),
\end{array}
$$
\item
$q=8$ and $f(x)g(x)=(x^2+4x+8)(x^4-16x^2+64)$,
\item
$q=2$ and $g(x)=x^4-4x^2+4$,
\item $q=2$ and $f(x)g(x)$ is one of the following polynomials
$$
(x^2+\eps2x+2)(x^4-\eps2x^2+4),\quad (x^2+2x+2)(x^4-2x^3+2x^2-2x+4).
$$
\end{itemize}
\end{teor}

\begin{teor}\label{teorfinal}
Let $f(x)$ be the product of three Frobenius polynomials of supersingular elliptic curves over $k$. If $q>4$ the isogeny class of abelian threefolds with Weil polynomial $f(x)$ contains a Jacobian, except in the following cases :
\begin{itemize}
\item $q=64$ and $f(x)=(x^2-16x+64)^3$,
\item $q=16$ and $f(x)$ is one of the following polynomials:
$$\ars{1.6}
\begin{array}{l}
(x^2+16)^2(x^2-8x+16),\, (x^2+16)(x^2\pm8x+16)^2,\, (x^2\pm8x+16)^3,\\
(x^2-8x+16)(x^2+8x+16)^2,\ (x^2-4x+16)(x^2\pm8x+16)^2,
\end{array}
$$
\item $q=8$ and $f(x)=(x^2-4x+8)^3$.
\end{itemize}

If $q=2$ these split isogeny classes never contain Jacobians. If $q=4$ the isogeny class of abelian threefolds with Weil polynomial $f(x)$ contains a Jacobian if and only if $f(x)$ is divisible by $(x^2+2x+4)(x^2-2x+4)$ or $f(x)$ is one of the following polynomials:
$$\ars{1.6}
\begin{array}{l}
(x^2+4)(x^2\pm2x+4)^2,\ (x^2+4x+4)(x^2\pm2x+4)^2,\ (x^2+4)^2(x^2+2x+4),\\
(x^2-4x+4)(x^2+4x+4)(x^2-2x+4).
\end{array}
$$
\end{teor}

\begin{cor}
If $q>64$ and the $k$-isogeny class of supersingular abelian threefolds with Weil polynomial $f(x)$ contains a Jacobian, then the $k$-isogeny class associated to the polynomial $f(-x)$ contains a Jacobian too.
\end{cor}

One may speculate if a similar result is valid for non-supersingular threefolds and/or other characteristics.

\begin{rem}
 Given an isogeny class which
contains a Jacobian, it is easy to use the constructive methods of
the proofs of Theorems \ref{teorsimple}, \ref{teorhalfsplit}, \ref{teorsplit} to build up explicit supersingular
quartics whose Jacobian lies in the prescribed isogeny class.
\end{rem}

If $q$ is a square, a curve is maximal (respectively minimal) if and
only if the Jacobian is $k$-isogenous to $E\times E\times E$ for $E$
an elliptic curve with Frobenius polynomial $x^2+2\sqrt{q}x+q$
(respectively $x^2-2\sqrt{q}x+q$).  Therefore, Theorem
\ref{teorfinal} solves in particular the problem of the existence of
maximal or minimal curves of genus three for $q$ square.

\begin{cor}
Let $q$ be a square. For $q>64$,  there are maximal curves and minimal curves of genus three over $k$. For $q=64$ there are maximal curves over $k$, but no minimal curve. For $q\le 16$ there are neither maximal nor minimal curves over $k$.
\end{cor}

\end{document}